\documentclass[12pt]{amsart}
\usepackage{amsfonts, amsbsy, amsmath, amssymb, cool}

\newtheorem{thm}{Theorem}[section]
\newtheorem{lem}[thm]{Lemma}

\newtheorem{rmk}[thm]{Remark}

\newtheorem{thm-con}[thm]{Theorem-Conjecture}
\numberwithin{equation}{section}

\theoremstyle{definition}

\newcommand{\f}{\Bbb F}

\begin{document}

\title[A note on Dickson polynomials of the third kind and  Legendre functions]{A note on Dickson polynomials of the third kind and  Legendre functions}

\author[Neranga Fernando]{Neranga Fernando}
\address{Department of Mathematics,
Northeastern University, Boston, MA 02115}
\email{w.fernando@northeastern.edu}

\author[Solomon Manukure]{Solomon Manukure}
\address{Department of Mathematics,
The University of Texas at Austin, Austin, TX 78712}
\email{smanukure@math.utexas.edu}

\begin{abstract}
In this paper, we show that the Dickson polynomials of the third kind satisfy a nonhomogeneous second order linear ordinary differential equation whose general solution contains Legendre functions.
\end{abstract}



\maketitle

\section{Introduction}

Let $n$ be a non-negative integer. It is well known that the two elementary symmetric polynomials $x_1+x_2$ and $x_1x_2$ form a $\mathbb{Z}$-basis of the ring of symmetric polynomials in $\mathbb{Z}[x_1,x_2]$. In \cite{Macdonald-1998}, it was shown that there exists a polynomial $D_n(x,y)\in \mathbb{Z}[x,y]$ such that 
$$x_1^n+x_2^n=D_n(x_1+x_2,x_1x_2).$$

The explicit form of $D_n(x,y)$ is given by the Waring's formula \cite[Theorem~1.1]{Lidl-Mullen-Turnwald-1993}

\[
D_{n}(x,y) = \sum_{i=0}^{\lfloor\frac n2\rfloor}\frac{n}{n-i}\dbinom{n-i}{i}(-y)^{i}x^{n-2i}.
\]

The polynomial $D_n(x,a)\in\,R[x]$ is the $n$-th Dickson polynomial of the first kind, where $R$ is a commutative ring with identity and $a\in R$ is a parameter. In \cite{Lidl-Mullen-Turnwald-1993}, it is shown that the Dickson polynomials of the first kind, $D_n(x,a)$, satisfy the homogeneous second order ordinary differential equaion
$$(x^2-4a)D_n^{\prime \prime}(x,a) + xD_n^{\prime}(x,a) - n^2\,D_n(x,a)=0.$$

The $n$-th Dickson polynomial of the second kind $E_n(x,a)$ is defined by
\[
E_{n}(x,a) = \sum_{i=0}^{\lfloor\frac n2\rfloor}\dbinom{n-i}{i}(-a)^{i}x^{n-2i},
\]
where $a\in R$ is a parameter. In \cite{Lidl-Mullen-Turnwald-1993}, it is also shown that the Dickson polynomials of the second kind, $E_n(x,a)$, satisfy the homogeneous second order ordinary differential equaion
$$(x^2-4a)E_n^{\prime \prime}(x,a) + 3xE_n^{\prime}(x,a) - n(n+2)\,E_n(x,a)=0.$$

Dickson polynomials have been widely studied over finite fields for their permutation behaviour. Let $\f_q$ denote the finite field with $q$ elements, where $q$ is a prime power. A permutation polynomial (PP) of $\Bbb F_q$ is a polynomial $f \in \Bbb F_q[x]$ such that the mapping $x\mapsto f(x)$ is a permutation of $\Bbb F_q$. Permutation polynomials over finite fields have important applications in coding theory, cryptography, finite geometry, combinatorics and computer science, among other fields. The reader may find an excellent presentation of Dickson polynomials in \cite{Lidl-Mullen-Turnwald-1993}. Even though most studies on PPs have been over finite fields, several authors have studied PPs over finite commutative rings. We refer the reader to \cite{GHM-FFA-2018}, \cite{Rivest-2001} and references therein for further details on PPs over finite commutative rings.

In \cite{Stoll-2008}, Stoll considered Dickson-type polynomials $f_n$ over $\mathbb{R}$, which depend on two real parameters $a$ and $B$. The explicit expression for $f_n$ is given by 

\begin{equation}\label{E1}
\displaystyle{f_n(x)=\sum_{i=0}^{\lfloor\frac n2\rfloor}\frac{n+(B-2)i}{n-i} \binom{n-i}{i}\,(-a)^{i}\,x^{n-2i}}.
\end{equation}

We point out that \eqref{E1} had also appeared earlier in \cite{Dujella-Gusic-2007}. 

Again in \cite{Stoll-2008}, Stoll showed that the Dickson-type polynomials $f_n$ with $a\neq 0, B\in \mathbb{R}$ satisfy the second order homogeneous ordinary differential equation (See \cite[Lemma 17]{Stoll-2008})

\begin{equation}\label{E2}
(A_4x^4+aA_2x^2+a^2A_0)f_n^{\prime \prime}+ (B_3x^3+aB_1x)f_n^{\prime}-(C_2x^2+aC_0)f_n=0,
\end{equation}

where $A_4, A_2, A_0, B_3, B_1, C_2, C_0 \in \mathbb{R}$. 

For $a\in \f_q$, the $n$-th Dickson polynomial of the $(k+1)$-th kind $D_{n,k}(x,a)$ as defined by Wang and Yucas in  \cite{Wang-Yucas-FFA-2012} is 
\begin{equation}\label{E11}
D_{n,k}(x,a) = \sum_{i=0}^{\lfloor\frac n2\rfloor}\frac{n-ki}{n-i}\dbinom{n-i}{i}(-a)^{i}x^{n-2i}.
\end{equation}

Note that all the coefficients of $D_{n,k}(x,a)$ are integers and $D_{n,k}(x,a)$ also satisfies \eqref{E2} with $A_4, A_2, A_0, B_3, B_1, C_2, C_0 \in \mathbb{Z}$. 

The associated Legendre equation \cite{Wysin-2011} is given by 

\begin{center}
$\displaystyle{(x^2-1)\,\frac{d^2u}{dx^2}+2x\,\frac{du}{dx}-\Big[l(l+1)+\frac{m^2}{x^2-1}\Big]\,u=0},$
\end{center}

where $l$ and $m$ are complex numbers. Its solutions are $m$-th derivatives of Legendre polynomials, which are a system of complete and orthogonal polynomials, multiplied by the adjustment factor $(1-x^2)^{\frac{m}{2}}$. These are the associated Legendre functions. The complex numbers $l$ and $m$ are called the degree and order of the associated Legendre functions, respectively. The Legendre polynomials are the associated Legendre functions of order $m=0$. We refer the reader to \cite{Arfken-Weber}, \cite{Oliver-Lozier-Boisvert-Clark}, \cite{Virchenko-Fedotova}, and \cite{Wysin-2011} for further details about Legendre functions.

Dickson polynomials have connections with special functions such as Legendre functions and hypergeometric functions. As stated previously, Dickson polynomials of the first and second kind satisfy some linear differential equations. It is well known that the solutions to these differential equations can be represented by Legendre functions. So far, there has been no work on how Dickson polynomials of the ($k+1$)-th kind are related to Legendre functions. In this note, our main goal is to explore the relationship between the Dickson polynomials of the third kind and associated Legendre functions. It turns out that the Dickson polynomials of the third kind satisfy a non-homogeneous second order linear ordinary differential equation whose solution contains Legendre functions. First, we show that the particular solution to the resulting differential equation involves the Dickson polynomials of the first kind and subsequently show that the general solution to the associated homogeneous equation involves the Associated Legendre functions which can be expressed in terms of the gamma function and hypergeometric functions.  

The paper is organized as follows. In Section 2, we show that the Dickson polynomials of the third kind satisfy the non-homogeneous second-order linear ordinary differential equation given in Lemma~\ref{L1.1}. In Section 3, we find its general solution by proving Theorem~\ref{T1.2}. 

By \eqref{E11}, the $n$-th Dickson polynomial of the third kind $D_{n,2}(x,a)$ is given by 

\begin{equation}\label{E1.4}
D_{n,2}(x,a) = \sum_{i=0}^{\lfloor\frac n2\rfloor}\frac{n-2i}{n-i}\dbinom{n-i}{i}(-a)^{i}x^{n-2i},
\end{equation}

where $a\in R$. 

Throughout the paper, we denote the $n$-th Dickson polynomial of the third kind $D_{n,2}(x,a)$ by $F_{n}(x,a)$. Also, we assume that $R=\mathbb{C}$. 

\begin{lem}\label{L1.1}

Let $a\ne 0$, $x=u+au^{-1}$ with $u\neq 0$ and $u^2\neq a$. Then we have

\begin{equation}\label{E1.1}
(x^2-4a)\,\,F_n^{\prime \prime}(x,a)+ 3x\,\,F_n^{\prime}(x,a)-n^2\,\,F_n(x,a)
= 2n\,\,D_n(x,a),
\end{equation}

\end{lem}

where $D_n(x,a)$ is the $n$-th Dickson polynomial of the first kind.

\vskip .2cm

In Section 3, we solve \eqref{E1.1} and prove the following theorem.

\begin{thm}\label{T1.2}
Let $F_c(x,a)$ and $F_p(x,a)$ be the general solution to the associated homogeneous equation of \eqref{E1.1} and the particular solution of \eqref{E1.1}, respectively. Then 

\[
\begin{split}
F_c(x,a)&=\frac{1}{\sqrt[4]{x^2-4a}}\,\Big[A\,P_{\sqrt{n^2+1}-\frac{1}{2}}^{(\frac{1}{2})}\Big(\frac{x}{2\sqrt{a}}\Big)+B\,Q_{\sqrt{n^2+1}-\frac{1}{2}}^{(\frac{1}{2})}\Big(\frac{x}{2\sqrt{a}}\Big)\Big],
\end{split}
\]

where $A$ and $B$ are constants, and $P_{\sqrt{n^2+1}-\frac{1}{2}}^{(\frac{1}{2})}$ and $Q_{\sqrt{n^2+1}-\frac{1}{2}}^{(\frac{1}{2})}$ are the associated Legendre functions of the first and second kinds, respectively. Also, 

$$F_p(x,a)=\displaystyle\sum_{k=0}^{n}\,b_k\,x^k,$$

where 

$$
\left\{
        \begin{array}{ll}
        b_{k+2}=-\Big(\displaystyle\frac{n^2-k\,(k+2)}{4\,a\,(k+2)\,(k+1)} \Big)\,b_k,\,\,\,\, n\,\, \text{is odd}, k\,\,\text{ is even}\,\,\, or \,\,\,n\,\, \text{is even}\,, k\,\,\text{ is odd} ,\\[0.5cm]
        b_{k+2}=\displaystyle\frac{1}{4\,a\,(k+2)\,(k+1)}\,\Big\{[k\,(k+2)-n^2]\,b_k-2n\,(\textnormal{coefficient of}\, x^{n-2i}\,\textnormal{in} \,D_n(x,a))\Big\},\\ otherwise,
        \end{array}
    \right.
$$
where $i=\frac{n-k}{2}$ and $D_n(x,a)$ is the $n$-th Dickson polynomial of the first kind.

In particular, we have 
$$F_n(x,a)=F_c(x,a)+F_p(x,a).$$

\end{thm}

\subsection{Definitions}

The Dickson polynomial $D_n(x,a)$ of the first kind of degree $n$ in the indeterminate $x$ and with parameter $a\in R$ is given as 

\[
D_{n}(x,a) = \sum_{i=0}^{\lfloor\frac n2\rfloor}\frac{n}{n-i}\dbinom{n-i}{i}(-a)^{i}x^{n-2i}. 
\]

Here $\lfloor\frac n2\rfloor$ denotes the largest integer $\leq n/2$. Let $u_1$ and $u_2$ be indeterminates. Then Waring's formula (see \cite[Theorem~1.1]{Lidl-Mullen-Turnwald-1993}) yields

\[
u_1^n+u_2^n = \sum_{i=0}^{\lfloor\frac n2\rfloor}\frac{n}{n-i}\dbinom{n-i}{i}(-u_1u_2)^{i}(u_1+u_2)^{n-2i}
\]

and thus

$$u_1^n+u_2^n=D_n(u_1+u_2,u_1u_2).$$

If we let $u_1=u$, $u_2=\frac{a}{u}$ and $x=u+\frac{a}{u}$, then the functional equation of the Dickson polynomial of the first kind is given by 

\begin{equation}\label{E1.6}
D_n(u+\frac{a}{u}, a)=u^{n}+\Big(\displaystyle\frac{a}{u}\Big)^n.
\end{equation}

For $a\ne 0$, Let $x=u+au^{-1}$ with $u\neq 0$ and $u^2\neq a$. Then the functional equation of the Dickson polynomial of the third kind, $F_n(x,a)$, is given by 

\begin{equation}\label{E1.5}
F_n(x,a)=(u+\frac{a}{u})\,\,\displaystyle\frac{u^n-\Big(\displaystyle\frac{a}{u}\Big)^n}{u-\displaystyle\frac{a}{u}},
\end{equation}

see \cite[Eq.~4.1]{Wang-Yucas-FFA-2012}. 

For $m$ a real number, the Leibniz rule for the product of two functions $u(z)$ and $v(z)$ is given by 

\begin{equation}\label{E1.8}
D_z^{m}uv\,=\,\displaystyle\sum_{\alpha=0}^{\infty}\,\binom{m}{\alpha}\,D_{z}^{m -\alpha}u\,D_z^{\alpha}v,
\end{equation}

where 

$$\binom{m}{\alpha}\,=\,\displaystyle\frac{\Gamma (m +1)}{\Gamma (m-\alpha +1)\,\Gamma (\alpha +1)}.$$


\section{Proof of Lemma~\ref{L1.1}}

From \eqref{E1.5}, we have

$$ F_n(u+au^{-1},a)\,\Big(u-\displaystyle\frac{a}{u}\Big)=(u+\frac{a}{u})\,\,\Big[u^n-\Big(\displaystyle\frac{a}{u}\Big)^n\Big].$$

By differentiating both sides with respect to $u$, we get

\[
\begin{split}
& F_n^{\prime}(u+au^{-1},a)\,\Big(u-\displaystyle\frac{a}{u}\Big)\,\,\Big(1-\displaystyle\frac{a}{u^2}\Big)+ F_n(u+au^{-1},a)\,\,\Big(1+\displaystyle\frac{a}{u^2}\Big)\cr
&= \Big(u+\displaystyle\frac{a}{u}\Big)\,\,\Big[nu^{n-1}+\displaystyle\frac{na^n}{u^{n+1}}\Big]+\Big[u^n-\Big(\displaystyle\frac{a}{u}\Big)^n\Big]\,\,\Big(1-\displaystyle\frac{a}{u^2}\Big).
\end{split}
\]

Multiplying by $u$ and differentiating both sides with respect to $u$ again yields

\[
\begin{split}
& F_n^{\prime \prime}(u+au^{-1},a)\,\Big(u-\displaystyle\frac{a}{u}\Big)^2\,\,\Big(1-\displaystyle\frac{a}{u^2}\Big)+ 2\,\,F_n^{\prime}(u+au^{-1},a)\,\,\Big(u-\displaystyle\frac{a}{u}\Big)\,\,\Big(1+\displaystyle\frac{a}{u^2}\Big)\cr
&+ F_n^{\prime}(u+au^{-1},a)\,\,\Big(u+\displaystyle\frac{a}{u}\Big)\,\,\Big(1-\displaystyle\frac{a}{u^2}\Big)+ F_n(u+au^{-1},a)\,\,\Big(1-\displaystyle\frac{a}{u^2}\Big) \cr
&= n\,\,\Big(u+\displaystyle\frac{a}{u}\Big)\,\,\Big[nu^{n-1}-\displaystyle\frac{na^n}{u^{n+1}}\Big]+n\,\,\Big[u^n+\Big(\displaystyle\frac{a}{u}\Big)^n\Big]\,\,\Big(1-\displaystyle\frac{a}{u^2}\Big)\cr
&+\Big[u^n-\Big(\displaystyle\frac{a}{u}\Big)^n\Big]\,\,\Big(1+\displaystyle\frac{a}{u^2}\Big)+\Big(u-\displaystyle\frac{a}{u}\Big)\,\,\Big[nu^{n-1}+\displaystyle\frac{na^n}{u^{n+1}}\Big].
\end{split}
\]

Multiplying by $u$ gives

\[
\begin{split}
& \Big(u-\displaystyle\frac{a}{u}\Big)^3\,\,F_n^{\prime \prime}(u+au^{-1},a)+ 3\,\,\Big[u^2-\Big(\displaystyle\frac{a}{u}\Big)^2\Big]\,\,F_n^{\prime}(u+au^{-1},a)\cr 
&+\Big(u-\displaystyle\frac{a}{u}\Big) \,\, F_n(u+au^{-1},a)\cr
&= (n^2+1)\,\,\Big(u-\displaystyle\frac{a}{u}\Big) \,\, F_n(u+au^{-1},a)+2n\,\,\Big(u-\displaystyle\frac{a}{u}\Big)\,\,\Big[u^{n}+\Big(\displaystyle\frac{a}{u}\Big)^n\Big], 
\end{split}
\]

which simplifies to 

\begin{equation}\label{E2.2}
\begin{split}
&\Big(u-\displaystyle\frac{a}{u}\Big)^2\,\,F_n^{\prime \prime}(u+au^{-1},a)+ 3\,\,\Big(u+\displaystyle\frac{a}{u}\Big)\,\,F_n^{\prime}(u+au^{-1},a)-n^2\,\,F_n(u+au^{-1},a)\cr
&= 2n\,\,\Big[u^{n}+\Big(\displaystyle\frac{a}{u}\Big)^n\Big]. 
\end{split}
\end{equation}

Note that, $D_n(x,a)=u^{n}+\Big(\displaystyle\frac{a}{u}\Big)^n$, where $D_n(x,a)$ is the $n$-th Dickson polynomial of the first kind (See \eqref{E1.6}). Also note that, $\Big(u-\displaystyle\frac{a}{u}\Big)^2=x^2-4a$. 

From \ref{E2.2} we have, 

$$
(x^2-4a)\,\,F_n^{\prime \prime}(x,a)+ 3x\,\,F_n^{\prime}(x,a)-n^2\,\,F_n(x,a)
= 2n\,\,D_n(x,a). 
$$

\begin{rmk}
Let $a\ne 0$ and $x=u+au^{-1}$ with $u\neq 0$. When $u^2= a$, the functional equation of $F_n(a,x)$ is given by 
\begin{equation}\label{E2.3}
F_n(x,a)=2(\pm\sqrt{a})^n\,n;
\end{equation}
see \cite[Remark 2.4]{Wang-Yucas-FFA-2012}. In this case, we have $F_n^{\prime}(x,a)=0$.
\end{rmk}

\section{Proof of Theorem~\ref{T1.2}}

In Section 2, we showed that Dickson polynomials of the third kind satisfy the second order non-homogeneous differential equation

\begin{equation}\label{E3.1}
(x^2-4a)\,\,F_n^{\prime \prime}(x,a)+ 3x\,\,F_n^{\prime}(x,a)-n^2\,\,F_n(x,a)
= 2n\,\,D_n(x,a). 
\end{equation}

We first find the particular solution to \eqref{E3.1}, $F_p(x,a)$. 

Let $F_p=\displaystyle\sum_{k=0}^{n}\,b_k\,x^k$ be the trial solution of the particular solution to \eqref{E3.1}. Then we have 
$$ F_p^{\prime}=\displaystyle\sum_{k=1}^{n}\,k\,b_k\,x^{k-1} \hspace{1cm}\textnormal{and}\hspace{1cm} F_p^{\prime \prime}=\displaystyle\sum_{k=2}^{n}\,k(k-1)\,b_k\,x^{k-2}.$$

Then 

\begin{equation}\label{E3.3}
\begin{split}
&(x^2-4a)\,\,F_n^{\prime \prime}(x,a)+ 3x\,\,F_n^{\prime}(x,a)-n^2\,\,F_n(x,a)\cr
&=(x^2-4a)\,\,\displaystyle\sum_{k=2}^{n}\,k(k-1)\,b_k\,x^{k-2}+ 3x\,\,\displaystyle\sum_{k=1}^{n}\,k\,b_k\,x^{k-1}-n^2\,\,\displaystyle\sum_{k=0}^{n}\,b_k\,x^k\cr
&=\displaystyle\sum_{k=0}^{n-2}\,\,[k(k-1)\,b_k+3\,k\,b_k-n^2\,b_k-4a(k+2)(k+1)\,b_{k+2}]\,x^k-b_{n-1}x^{n-1}+2nb_nx^n
\end{split}
\end{equation}

Let $b_{n+1}=b_{n+2}=0.$ Now \eqref{E3.3} can be written as

\[
\begin{split}
&(x^2-4a)\,\,F_n^{\prime \prime}(x,a)+ 3x\,\,F_n^{\prime}(x,a)-n^2\,\,F_n(x,a)\cr
&=\displaystyle\sum_{k=0}^{n}\,\,[k(k-1)\,b_k+3\,k\,b_k-n^2\,b_k-4a(k+2)(k+1)\,b_{k+2}]\,x^k
\end{split}
\]

From \eqref{E3.1} we have

\begin{equation}\label{E3.4}
\begin{split}
&\displaystyle\sum_{k=0}^{n}\,\,[k(k-1)\,b_k+3\,k\,b_k-n^2\,b_k-4a(k+2)(k+1)\,b_{k+2}]\,x^k\cr 
&= 2n\,\sum_{i=0}^{\lfloor\frac n2\rfloor}\frac{n}{n-i}\dbinom{n-i}{i}(-a)^{i}x^{n-2i}.
\end{split}
\end{equation}

Consider the following cases. 

\begin{itemize}
\item [(i)] $n$ is odd and $k$ is even, 
\item [(ii)] $n$ is even and $k$ is odd, 
\item [(iii)]  $n$ is odd and $k$ is odd, and
\item [(iv)]  $n$ is even and $k$ is even.
\end{itemize}

In cases (i) and (ii), by comparing the coefficients in \eqref{E3.4} we have
$$k(k-1)\,b_k+3\,k\,b_k-n^2\,b_k-4a(k+2)(k+1)\,b_{k+2}=0$$
which implies 
$$b_{k+2}=-\Big(\displaystyle\frac{n^2-k\,(k+2)}{4\,a\,(k+2)\,(k+1)} \Big)\,b_k.$$

Now consider cases (iii) and (iv). By comparing the coefficients in \eqref{E3.4} we have

$$b_{k+2}=\displaystyle\frac{1}{4\,a\,(k+2)\,(k+1)}\,\Big\{[k\,(k+2)-n^2]\,b_k-2n\,(\textnormal{coefficient of}\, x^{n-2i}\,\textnormal{in} \,D_n(x,a))\Big\},$$

where $i=\frac{n-k}{2}$. 

Now let's consider the associated homogeneous equation, i.e.

\begin{equation}\label{E3.2}
(x^2-4a)\,\,F_n^{\prime \prime}(x,a)+ 3x\,\,F_n^{\prime}(x,a)-n^2\,\,F_n(x,a)= 0. \\[.5cm]
\end{equation}

\noindent Denote $F_n(x,a)$ by $F(x)$ and let $X=\displaystyle\frac{x}{2\sqrt{a}}$. It then follows that \\

$$\displaystyle{2\sqrt{a}\frac{d}{dX}\Big[\Big(X^2-1\Big)^{\frac{3}{2}}\,F^{\prime}(2\sqrt{a}X)\Big]-n^2\,\Big(X^2-1\Big)^{\frac{1}{2}}\,F(2\sqrt{a}X)=0}.$$

\noindent Letting $P(X)=F(2\sqrt{a}X)$ gives

\begin{equation}\label{E3.5}
\displaystyle{\frac{d}{dX}\Big[\Big(X^2-1\Big)^{\frac{3}{2}}\,P^{\prime}(X)\Big]-n^2\,\Big(X^2-1\Big)^{\frac{1}{2}}\,P(X)=0}.
\end{equation}

Now consider the Legendre equation:

$$(1-x^2)P_l^{\prime \prime}(x)-2x\,P_l^{\prime}(x)+l(l+1)P_l(x)=0$$

or equivalently

\begin{equation}\label{E3.8}
(x^2-1)P_l^{\prime \prime}(x)+2x\,P_l^{\prime}(x)-l(l+1)P_l(x)=0,
\end{equation}

where $l\in \mathbb{R}$ with $l\geq 0$.

Differentiating $m(\in \mathbb{R})$ times using Leibniz rule (see \eqref{E1.8}) gives the following. 

\begin{equation}\label{E3.6}
\begin{split}
&D_x^m[(x^2-1)\,P_l^{\prime \prime}(x)]=\sum\limits_{\alpha=0}^{\infty}\frac{\Gamma(m+1)}{\Gamma(m-\alpha+1)\Gamma(\alpha+1)}D_x^{\alpha}(x^2-1)D^{m-\alpha}P_l^{\prime \prime}(x)\cr
&=(x^2-1)P_l^{(m+2)}(x)+2mx\,\,P_l^{(m+1)}(x)+m(m-1)\,\,P_l^{(m)}(x).
\end{split}
\end{equation}

Similarly, 

\begin{equation}\label{E3.7}
\displaystyle{D_x^m[2x\,P_l^{\prime}(x)]=2x\,P_l^{(m+1)}(x)+2m\,P_l^{(m)}(x)}.
\end{equation}

\eqref{E3.8}, \eqref{E3.6}, and \eqref{E3.7} yield

\begin{equation}\label{E3.9}
(x^2-1)P_l^{(m+2)}(x)+2(m+1)\,x\,\,P_l^{(m+1)}(x)-(l-m)(l+m+1)P_l^{(m)}(x)=0.
\end{equation}

Now let $y(x)=P_l^{(m)}$. Then \eqref{E3.9} becomes

\begin{equation}\label{E3.10}
(x^2-1)\frac{d^2y}{dx^2}+2(m+1)\,x\,\,\frac{dy}{dx}-(l-m)(l+m+1)\,y=0.
\end{equation}

By multiplying \eqref{E3.10} by $(x^2-1)^m$, we get

\begin{equation}\label{E3.11}
\displaystyle{\frac{d}{dx}\Big[\Big(x^2-1\Big)^{m+1}\,\frac{dy}{dx}\Big]-(l-m)(l+m+1)\,\Big(x^2-1\Big)^{m}\,y=0}.
\end{equation}

Now let 

\begin{equation}\label{Solo}
y(x)=(x^2-1)^{-\frac{m}{2}}\,u(x). 
\end{equation}

Then $$\displaystyle\frac{dy}{dx}=(x^2-1)^{-\frac{m}{2}}\,\frac{du}{dx}-mx(x^2-1)^{-\frac{m}{2}-1}\,u.$$

Thus \eqref{E3.11} becomes

\[
\begin{split}
&\Big(x^2-1\Big)^{\frac{m}{2}+1}\,\frac{d^2u}{dx^2}+\Big(\frac{m}{2}+1\Big)(2x)\Big(x^2-1\Big)^{\frac{m}{2}}\,\frac{du}{dx}-mx\Big(x^2-1\Big)^{\frac{m}{2}}\,\frac{du}{dx}-m\Big(x^2-1\Big)^{\frac{m}{2}}\,u\cr
&-mx\,\frac{m}{2}\,(2x)\,\Big(x^2-1\big)^{\frac{m}{2}-1}\,u-(l-m)(l+m+1)\,\Big(x^2-1\Big)^{\frac{m}{2}}\,u=0. 
\end{split}
\]

Dividing by $\Big(x^2-1\Big)^{\frac{m}{2}}$, we obtain

\begin{equation}\label{E3.12}
\displaystyle{(x^2-1)\,\frac{d^2u}{dx^2}+2x\,\frac{du}{dx}-\Big[l(l+1)+\frac{m^2}{x^2-1}\Big]\,u=0},
\end{equation}

which is the associated Legendre differential equation, with $0\leq m\leq l$. Its general solution is given by (see \cite[Ch.~5]{Wang-Guo})
$$u_g=c_1\,P_l^{(m)}(x)\,+\,c_2\,Q_l^{(m)}(x),$$

where $c_1$ and $c_2$ are arbitrary constants, and $P_l^{(m)}(x)$ and $Q_l^{(m)}(x)$ are the associated Legendre functions of the first and second kinds, respectively. Thus, according to \eqref{Solo}
$$y_g(x)=\Big(x^2-1\Big)^{-\frac{m}{2}}\,u_g(x)$$
is the general solution to \eqref{E3.11}. 

Now consider the equation

\begin{equation}\label{E3.13}
\frac{d}{dx}\Big[\Big(x^2-1\Big)^{\frac{3}{2}}\,\frac{dy}{dx}\Big]-n^2\,\Big(x^2-1\Big)^{\frac{1}{2}}\,y=0.
\end{equation}

This is in the form of \eqref{E3.11} with $m=\frac{1}{2}$ and $(l-m)(l+m+1)=n^2$. Since $l\in \mathbb{R}$ with $l\geq 0$, we have $l=\sqrt{n^2+1}-\frac{1}{2}$. Then the solution to \eqref{E3.13} is 

$$y=(x^2-1)^{-\frac{1}{4}}\,\Big[c_1\,P_{\sqrt{n^2+1}-\frac{1}{2}}^{(\frac{1}{2})}(x)+c_2\,Q_{\sqrt{n^2+1}-\frac{1}{2}}^{(\frac{1}{2})}(x)\Big].$$

So the solution to  \eqref{E3.5} is 

$$P(X)=(X^2-1)^{-\frac{1}{4}}\,\Big[c_1\,P_{\sqrt{n^2+1}-\frac{1}{2}}^{(\frac{1}{2})}(X)+c_2\,Q_{\sqrt{n^2+1}-\frac{1}{2}}^{(\frac{1}{2})}(X)\Big].$$

Since  $X=\displaystyle\frac{x}{2\sqrt{a}}$ and $P(X)=F(2\sqrt{a}X)$, we have

\[
\begin{split}
F(x)&=(x^2-4a)^{-\frac{1}{4}}\,\Big[\frac{c_1}{(4a)^{-\frac{1}{4}}}\,P_{\sqrt{n^2+1}-\frac{1}{2}}^{(\frac{1}{2})}\Big(\frac{x}{2\sqrt{a}}\Big)+\frac{c_2}{(4a)^{-\frac{1}{4}}}\,Q_{\sqrt{n^2+1}-\frac{1}{2}}^{(\frac{1}{2})}\Big(\frac{x}{2\sqrt{a}}\Big)\Big].
\end{split}
\]

If we denote the general solution to \eqref{E3.2} by $F_c(x,a)$, then we have

\[
\begin{split}
F_c(x,a)&=\frac{1}{\sqrt[4]{x^2-4a}}\,\Big[A\,P_{\sqrt{n^2+1}-\frac{1}{2}}^{(\frac{1}{2})}\Big(\frac{x}{2\sqrt{a}}\Big)+B\,Q_{\sqrt{n^2+1}-\frac{1}{2}}^{(\frac{1}{2})}\Big(\frac{x}{2\sqrt{a}}\Big)\Big]
\end{split}
\]

where $A$ and $B$ are arbitrary constants. Hence the general solution to the differential equation

\[
(x^2-4a)\,\,F_n^{\prime \prime}(x,a)+ 3x\,\,F_n^{\prime}(x,a)-n^2\,\,F_n(x,a)
= 2n\,\,D_n(x,a),
\]

is

\begin{equation}\label{E4.3}
\begin{split}
F_n(x,a)&=\frac{1}{\sqrt[4]{x^2-4a}}\,\Big[A\,P_{\sqrt{n^2+1}-\frac{1}{2}}^{(\frac{1}{2})}\Big(\frac{x}{2\sqrt{a}}\Big)+B\,Q_{\sqrt{n^2+1}-\frac{1}{2}}^{(\frac{1}{2})}\Big(\frac{x}{2\sqrt{a}}\Big)\Big] + \displaystyle\sum_{k=0}^{n}\,b_k\,x^k, 
\end{split}
\end{equation}

where $b_k$ is determined as explained in Theorem~\ref{T1.2}. 

In terms of hypergeometric functions, the Legendre functions can be written as (see \cite[Ch.~5]{Wang-Guo})

\begin{equation}\label{E4.1}
\begin{split}
P_{\sqrt{n^2+1}-\frac{1}{2}}^{(\frac{1}{2})}(z)=\frac{1}{\Gamma (\frac{1}{2})}\,\,\Big[\frac{1+z}{1-z}\Big]^{\frac{1}{4}}\,\,\Hypergeometric{2}{1}{-\sqrt{n^2+1}+\frac{1}{2},\sqrt{n^2+1}+\frac{1}{2}}{\frac{1}{2}}{\frac{1-z}{2}},
\end{split}
\end{equation}

for $\mid 1-z \mid <2$, and 

\begin{equation}\label{E4.2}
\begin{split}
&Q_{\sqrt{n^2+1}-\frac{1}{2}}^{(\frac{1}{2})}(z)\cr
&=\frac{\sqrt{\pi}}{2^{\sqrt{n^2+1}+\frac{1}{2}}}\,\,\frac{i\,(z^2-1)^{\frac{1}{4}}}{z^{\sqrt{n^2+1}+1}}\,\,\Hypergeometric{2}{1}{\frac{\sqrt{n^2+1}+1}{2},\frac{\sqrt{n^2+1}+2}{2}}{\sqrt{n^2+1}+1}{\frac{1}{z^2}}
\end{split}
\end{equation}

for $\mid z \mid > 1$, 

where $\Gamma$ denotes the gamma function and $\Hypergeometric{2}{1}{a,b}{c}{z}$ is the hypergeometric function defined by the power series

$$\displaystyle{\Hypergeometric{2}{1}{a,b}{c}{z}= \sum_{n=0}^{\infty}\frac{(a)_n\,(b)_n}{(c)_n} \frac{z^n}{n!},}$$

for $\mid z \mid <1$.

Here $(\cdot)_n$ denotes the rising factorial defined by

$$
(a)_n = \left\{
        \begin{array}{ll}
            1 & ,\quad n=0 \\[0.2cm]
            a(a+1)\cdots (a+n-1)& ,\quad n>0.
        \end{array}
    \right.
$$

\section{Concluding remarks}
In this note we have shown that there is a relationship between Dickson polynomials of the third kind and Legendre functions through the existence of a second order linear differential equation. Currently we do not know if such a differential equation exists for the Dickson polynomials of the ($k+1$)-th kind. The nature of our computations clearly indicates that finding such a differential equation, if it exists, may be an arduous exercise.

\end{document}